\documentclass[11pt]{amsart2000}
\usepackage{amssymb}
\newtheorem{thm}{Theorem}[section]
\newtheorem{fac}[thm]{Fact}
\newtheorem{pro}[thm]{Proposition}
\newtheorem{lem}[thm]{Lemma}
\newtheorem{cor}[thm]{Corollary}

\newtheorem{claim}{Claim}[thm]

\begin{document}
\setlength{\unitlength}{0.01in}
\linethickness{0.01in}
\begin{center}
\begin{picture}(474,66)(0,0)
\multiput(0,66)(1,0){40}{\line(0,-1){24}}
\multiput(43,65)(1,-1){24}{\line(0,-1){40}}
\multiput(1,39)(1,-1){40}{\line(1,0){24}}
\multiput(70,2)(1,1){24}{\line(0,1){40}}
\multiput(72,0)(1,1){24}{\line(1,0){40}}
\multiput(97,66)(1,0){40}{\line(0,-1){40}}
\put(143,66){\makebox(0,0)[tl]{\footnotesize Proceedings of the Ninth Prague Topological Symposium}}
\put(143,50){\makebox(0,0)[tl]{\footnotesize Contributed papers from the symposium held in}}
\put(143,34){\makebox(0,0)[tl]{\footnotesize Prague, Czech Republic, August 19--25, 2001}}
\end{picture}
\end{center}
\vspace{0.25in}
\setcounter{page}{195}
\title[Nets in Ordered Fields]{Using nets in Dedekind, monotone, or Scott
incomplete ordered fields and definability issues}
\author{Mojtaba Moniri}
\author{Jafar S. Eivazloo}
\address{Department of Mathematics, Tarbiat Modarres University, Tehran,
Iran}
\email{moniri\_m@modares.ac.ir}
\email{eivazl\_j@modares.ac.ir}
\thanks{Mojtaba Moniri and Jafar S. Eivazloo,
{\em Using nets in Dedekind, monotone, or Scott incomplete ordered fields 
and definability issues},
Proceedings of the Ninth Prague Topological Symposium, (Prague, 2001),
pp.~195--203, Topology Atlas, Toronto, 2002}
\begin{abstract}
Given a Dedekind incomplete ordered field, a pair of convergent nets of
gaps which are respectively increasing or decreasing to the same point is
used to obtain a further equivalent criterion for Dedekind completeness of
ordered fields: Every continuous one-to-one function defined on a closed
bounded interval maps interior of that interval to the interior of the
image. 
Next, it is shown that over all closed bounded intervals in any monotone 
incomplete ordered field, there are continuous not uniformly continuous
unbounded functions whose ranges are not closed, and continuous 1-1
functions which map every interior point to an interior point (of the
image) but are not open. 
These are achieved using appropriate nets cofinal in gaps or coinitial in
their complements. 
In our third main theorem, an ordered field is constructed which has
parametrically definable regular gaps but no $\emptyset$-definable
divergent Cauchy functions (while we show that, in either of the two cases
where parameters are or are not allowed, any definable divergent Cauchy
function gives rise to a definable regular gap). 
Our proof for the mentioned independence result uses existence of infinite
primes in the subring of the ordered field of generalized power series
with rational exponents and real coefficients consisting of series with no
infinitesimal terms, as recently established by D.~Pitteloud.
\end{abstract}
\maketitle

\section{A Dedekind Incompleteness Feature via Convergent Nets of Gaps}

A cut of an ordered field $F$ is a subset which is downward closed in
$F$. By a nontrivial cut, we mean a nonempty proper cut. 
A nontrivial cut is a gap if it does not have a least upper bound in the
field. 
An ordered field is Archimedean (has no infinitesimals) just in case it
can be embedded in ${\mathbb R}$.

The following fact presents some of the well known characterizations of
the ordered field of real numbers. 
A more delicate equivalent condition is presented in Theorem 1.2.

\begin{fac}\label{easy2nd}
The real ordered field ${\mathbb R}$ is, up to isomorphism, the
unique ordered field which satisfies either of the following
equivalent conditions:
\begin{itemize}
\item[(i)] Dedekind Completeness, i.e. not having any gaps,
\item[(ii)] connectedness,
\item[(iii)] not being totally disconnected,
\item[(iv)] every (nonempty) convex subset being an interval (of one of
the usual kinds),
\item[(v)] all convex subsets being connected,
\item[(vi)] all intervals being connected,
\item[(vii)] all continuous functions on the field mapping any convex
subset onto a convex subset,
\item[(viii)] all continuous functions on the field satisfying the
intermediate value property.
\end{itemize}
\end{fac}

These are fairly well known, let us give the argument for (iii).
Pick a Dedekind incomplete ordered field $F$. In the Archimedean
case, $F$ is a proper subfield of ${\mathbb R}$ which therefore
misses some points in any real interval. In the non-Archimedean
case, and given any two points $a<b$ of a subset $A$ of $F$, we
have the nontrivial clopen subset of those points $x$ in $A$ such
that $\frac{x-a}{b-a}$ is a ${\mathbb Q}$-infinitesimal.
Alternatively, one can use linear increasing functions between
intervals to map a given gap (somewhere in the field) into a given
interval.

Notice that if $F$ is an ordered field and $P(x)\in F[x]$, then
$P(x)$ is continuous on $F$. To see this, it suffices to show that
if $a\in F$ and $P(a)<0$, then $P$ is negative on a neighborhood
of $a$ in $F$. As the same inequality also holds in $RC(F)$, where
polynomials have factorizations into linear and irreducible
quadratics, there are $b,c\in RC(F)$ with $b<a<c$ such that
$\forall x\in (b,c)_{RC(F)}$, $P(x)<0$. Using cofinality of $F$ in
$RC(F)$, one can now see that there are $d,e\in F$ such that
$b<d<a<e<c$ (observe that there exists $g\in F$ such that
$0<g<a-b$ and put $d=a-g$, similarly for $e$). A very mild use of
this takes place in the claim within the next theorem.

\begin{thm}\label{involved2nd}
An ordered field is Dedekind complete if and only if all
continuous 1-1 functions defined on some (equivalently all)
non-degenerated closed bounded interval(s) of it map interior
points [of the interval(s)] to interior points [of their
range(s)].
\end{thm}

\begin{proof}

We only need to prove the {\it if} part.

\begin{claim}
Given two convex subsets $A$ and $B$ of an ordered field $F$, with $A$
bounded and $B$ not a singleton, there exists a nonconstant linear (and so
1-1 continuous) function from $A$ into $B$.
\end{claim}

\begin{proof}[Proof of Claim] 
Let $A\subseteq [a_{1},a_{2}]$, $[b_{1},b_{2}]\subseteq B$, with
$a_{1}\not =a_{2}$ and $b_{1}\not =b_{2}$. 
Consider the function 
$$f(x)=b_{1}+\frac{b_{2}-b_{1}}{a_{2}-a_{1}}(x-a_{1}),$$ 
whose restriction to $A$ does the required job.
\end{proof}

Let $F$ be a Dedekind incomplete ordered field of cofinality $\lambda$. 
Consider an interval $[a,b]$ in $F$ with a point $c\in (a,b)$. 
Pick a strictly increasing net $(a_{\alpha})_{\alpha <\lambda}$ and a
strictly decreasing one $(b_{\alpha})_{\alpha <\lambda}$ in $[a,b]$ with
$a_{0}=a$ and $b_{0}=b$ which converge in $F$ to $c$. 
For each $\alpha$ less than $\lambda$, let $U_{\alpha}$ be a gap in
$(a_{\alpha },a_{\alpha +1})$ and $V_{\alpha}$ a gap in 
$(b_{\alpha +1 },b_{\alpha })$. 
Put 
$$
S_{0}=U_{0}\cap [a,b],\ 
T_{0}=[a,b]\setminus V_{0},$$
$$
S_{\alpha+1}=U_{\alpha+1}\setminus U_{\alpha},\
T_{\alpha+1}=V_{\alpha}\setminus V_{\alpha +1}$$ 
for $\alpha < \lambda$, and for limit $\beta<\lambda$, 
$$
S_{\beta}=U_{\beta}\setminus\cup_{\gamma<\beta}U_{\gamma},\
T_{\beta}=(\cap_{\gamma<\beta}V_{\gamma})\setminus V_{\beta}.
$$
Note that $[a,c)=\cup_{\alpha<\lambda}S_{\alpha}$, and 
$(c,b]=\cup_{\alpha<\lambda}T_{\alpha}$, and that both unions are 
disjoint. 
Using the claim, consider a function $f$ on $[a,b]$ which linearly
injects, for $\alpha <\lambda$, the set $S_{\alpha}$ into $T_{\alpha \cdot
2}$, $T_{\alpha}$ into $T_{(\alpha\cdot 2 )+1}$, and finally $c$ to $c$. 
Then $f$ maps the interior point $c$ of $[a,b]$ to the boundary point $c$
of $f([a,b])$. 
\end{proof}

\begin{cor}\label{openmaps}
An ordered field is Dedekind complete if and only if all continuous 1-1
functions on some (equivalently all) non-degenerated closed bounded
interval(s) are open.
\end{cor}

\begin{proof}
For the non-trivial {\em if} part, notice that all open maps on any
ordered field map interior to interior.
\end{proof}

\section{Some Consequences of Monotone Incompleteness via Nets Cofinal in
Gaps}

An ordered field $F$ is Scott complete, see \cite{S}, if it
satisfies either of the following three equivalent conditions:
\begin{itemize}
\item[(i)]
It does not have any proper extensions to an ordered field in
which it is dense,
\item[(ii)]
It does not have any regular gaps, where a gap $C\subset F$
is regular if $\forall\epsilon\in F^{>0}$,
$C+\epsilon\not\subseteq C$,
\item[(iii)]
All functions $f:F\rightarrow F$ which are Cauchy at
$\infty_{F}$, are convergent there.
\end{itemize}

We abbreviate the (equivalent non first order) axiomatizations obtained
from the theory $OF$ of ordered fields by adding condition (ii) or
(iii) as $S_{rc}COF$ and $S_{cf}COF$ respectively.

Observe that an ordered field of cofinality $\lambda$ is Scott complete if
and only if it has no divergent Cauchy nets of length $\lambda$.

If an ordered field $F$ is not Scott complete, then for all $a<b\in F$,
there is a regular gap $V$ of $F$ with $a\in V$ and $b\in F\setminus V$
(and so a divergent Cauchy net in $(a,b)$ of length equal to cofinality of
the field). In fact, any regular gap of $F$ has an additive translation to
a regular gap in $(a,b)$. To see this, let $U$ be a regular gap of $F$. We
can pick $c\in U$ and $d\in F\setminus U$ such that $d-c<b-a$. Let
$V=U+(a-c)$. It is straightforward to see that $V$ is a regular gap and
$a<V<b$.

It was proved in \cite{S} (Theorem 1), that any ordered field $F$ has a
(unique, up to an isomorphism of ordered fields which is identity on $F$)
Scott completion. It is characterized by being Scott complete and having
$F$ dense in it. Furthermore $F$ is dense in $RC(F)$ if and only if its
Scott completion is real closed, see \cite{S} (Theorem 2).

As shown in \cite{EGH} (Lemma 2.3), any uncountable ordered field has a
dense transcendence basis over the rationals. Therefore all Scott complete
ordered fields necessarily have proper dense subfields and so can be
obtained by Scott completion of a proper subfield (the field extension,
inside the original ordered field, of rationals by the just mentioned
basis minus a point does the job).

We will state and prove a number of results on Scott completeness and its
two first order versions in section 3. In this section, we are aiming at a
result dealing with the stronger notion of monotone completeness. Monotone
complete ordered fields were introduced in \cite{K}. They are ordered
fields with no bounded strictly increasing divergent functions. Scott
complete ordered fields that are $\kappa$-Archimedean for some regular
cardinal $\kappa$ are Monotone Complete (by $\kappa$-Archimedean, one
means that there are subsets of cardinality $\kappa$ whose distinct
elements are at least a unit apart and all such sets are unbounded), see
\cite{CL}.

As we have already mentioned, all open maps on any interval in an
arbitrary ordered field map every interior point to an interior point of
the image. However, the converse property, that of mapping interior to
interior implying being open even when restricted to continuous 1-1
functions, is strong enough to imply monotone completeness as the next
theorem shows.

\begin{thm} 
Over all closed bounded intervals in arbitrary monotone incomplete ordered
fields, there are
\begin{itemize}
\item[(i)]
continuous non uniformly continuous unbounded functions whose ranges are
not closed,
\item[(ii)]
continuous 1-1 functions mapping all interior points to interior of the
image which, nevertheless, are not open.
\end{itemize}
\end{thm}

\begin{proof}
Let $F$ be a monotone incomplete ordered field of cofinality $\lambda$. 
We first observe the following.

\begin{claim}
Any nondegenerated interval of $F$ contains a strictly increasing
divergent net of length $\lambda$ [compare with \cite{K} (Proposition
1(a))].
\end{claim}

\begin{proof}[Proof of Claim] 
Given a divergent bounded strictly increasing net of length $\lambda$,
consider an interval containing its range. 
Then for any other interval, one may apply the linear strictly increasing
function mapping the former onto the latter interval, thereby producing a
strictly increasing divergent net $\lambda$ in the latter. 
This concludes the claim.
\end{proof}

(i). 
Given a closed bounded interval $[a,b]$ of $F$, using the claim, we can
pick divergent $\lambda$-nets $(a_{\alpha})_{\alpha <\lambda}$ and
$(b_{\alpha})_{\alpha <\lambda}$ with $a_{0}=a$, $b_{0}=b$,
$a_{\alpha}<b_{\beta}$, for all $\alpha,\beta <\lambda$, which are
respectively strictly increasing or strictly decreasing. 
Also, let $(u_{\alpha})_{\alpha <\lambda}$ and $(d_{\alpha})_{\alpha
<\lambda}$ be both strictly increasing, the former being unbounded and
having consecutive terms which are at least a unit apart, the latter
converging in $F$ to $0$. 
Now, for each $\alpha<\lambda$, let $f$ map $[a_{\alpha},a_{\alpha +1})$
linearly and increasingly onto $[u_{\alpha},u_{\alpha +1})$ and 
$(b_{\alpha +1},b_{\alpha}]$ linearly and in a decreasing manner onto
$[d_{\alpha},d_{\alpha+1})$ and be equal to 1 on the rest of $[a,b]$. 
We may assume without loss of generality that for any $\delta>0$, there
exists $\alpha_{0} <\lambda$ such that $a_{\alpha_{0}
+1}-a_{\alpha_{0}}<\delta$. 
To see that such a reduction indeed causes no loss of generality, take a
net $(c_{\alpha})_{\alpha<\lambda}$ converging to $0$ such that
$0<c_{\alpha}<a_{\alpha +1} - a_{\alpha}$ and consider the net obtained
from $(a_{\alpha})_{\alpha<\lambda}$ by right-shifting those of its terms
which have limit ordinal indices by the corresponding same-indexed terms
in $(c_{\alpha})_{\alpha<\lambda}$. 
The function $f$ is a continuous non uniformly continuous function on
$[a,b]$ whose values are $F$-unbounded on the downward closure of
$(a_{\alpha})_{\alpha <\lambda}$ in $[a,b]$ and asymptotic to zero on the
upward closure of $(b_{\alpha})_{\alpha <\lambda}$ in $[a,b]$ traversed
backwards. 
To see that $f$ is not uniformly continuous, let $\epsilon =\frac{1}{2}$. 
Then for any $\delta >0$, by the assumption we made above, there exists
$\alpha_{0} <\lambda$ such that $a_{\alpha_{0} +1}-a_{\alpha_{0}}<\delta$. 
Now by $u_{\alpha +1}-u_{\alpha} \geq 1$ and the way $f$ is constructed,
we have 
$$f(a_{\alpha_{0} +1})-f(a_{\alpha_{0}})\geq 1 >\epsilon.$$

(ii).
Consider the interval $[a,b]$ with a point $c\in (a,b)$. 
We now follow a construction similar to the one in Theorem 
~\ref{involved2nd}. 
Using the same notation as there, the change we make is the following. 
We linearly inject $T_{\alpha}$, for $1\leq \alpha <\lambda$, into
$T_{(\alpha\cdot 2 )+1}$ but $T_{0}$ onto $[a,c)$. 
The latter can be done as in part (i) in a piece-wise manner. 
Here we may assume that $T_{0}$ is the upward closure of a strictly
decreasing divergent $\lambda$-net in $(b_{1},b_{0}]$. 
Then $c$ will be a boundary point of the image of the open interval
$(a_{1},b_{1})$ under $f$, so that image is not open. 
On the other hand, by construction, $f$ sends interior to interior.
\end{proof}

\section{P-Definable Regular Gaps Not Traversed by $\emptyset$-Definable
Cauchy Functions}
	
We now consider two notions of being definably (with or without
parameters) Scott complete for ordered fields, those corresponding to
$S_{rc}COF$ and $S_{cf}COF$. They are ordered fields with no definable
regular gaps, respectively no definable functions which are Cauchy at
infinity but divergent there. We denote the corresponding theories by
$D_{p}S_{rc}COF$, $D_{p}S_{cf}COF$, $D_{\emptyset}S_{rc}COF$, and
$D_{\emptyset}S_{cf}COF$.

It is easy to see that, as long as there are no definability concerns, all
regular gaps can be traversed by suitable (divergent Cauchy) functions and
vice versa, every divergent Cauchy function induces a regular gap. The
latter converse is shown in Theorem ~\ref{conditional}(i) to be true in
either of the p-definable or $\emptyset$-definable cases. For the former
direction, however, we are only able in Theorem ~\ref{conditional}(ii) to
prove an independence result in a mixed $\emptyset$-definable /
p-definable case.

\begin{lem}\label{shepherdson}
If the ordered field $F$ is a proper dense sub-field of its real closure,
then $F\nvDash D_{p}S_{rc}COF$.
\end{lem}

\begin{proof}
Pick out an element $r\in RC(F)\setminus F$ and consider the set 
$$C=\{x\in F: \ RC(F)\vDash x<r\}.$$
It is obviously a gap of $F$. As $F$ is dense in $RC(F)$, $C$ is regular. 
Finally, $C$ is p-definable in $F$: consider the minimal polynomial of $r$
over $F$ and the number of roots of that polynomial in $RC(F)$ which are
less than $r$.
\end{proof}

For any ordered field $F$ and ordered abelian group $G$, the set
$[[F^{G}]]$ of all functions $G\rightarrow F$ whose supports are well
ordered in $G$ equipped with pointwise sum and Cauchy product
$$(f_{1}f_{2})(g)=\sum_{i+j=g}f_{1}(i)f_{2}(j)$$ 
(a finite sum by the condition on the supports) forms a field. It can be 
ordered by comparison of values at the minimum of support of the 
difference. 
Elements of $[[F^{G}]]$ can also be thought of as those formal power 
series $\sum_{g\in G}f(g)t^{g}$ which have well ordered supports. 
The indeterminate $t$ is taken to be a positive $F$-infinitesimal. The 
field $[[F^{G}]]$ is real closed if and only if $F$ is so and $G$ is 
divisible, see \cite{R} (6.10).

In the proposition below, $\chi_{S}$ stands for the characteristic 
function of (a set) $S$.

\begin{pro}\label{scottfieldsofgps}
For any ordered abelian group $G$ and ordered field $F$, the generalized
power series field $[[F^{G}]]$ is Scott complete if $F$ is so.
\end{pro}

\begin{proof}
Let $\lambda$ be the cofinality of $G$ and consider a strictly decreasing
$\lambda$-net $(a_{\theta})_{\theta <\lambda}$ in $G^{<0}$ which is
unbounded below there. 
Assume $\mathcal{F}:[[F^{G}]]\rightarrow [[F^{G}]]$ is Cauchy. 
For each $g\in G$, define $f_{g}:\lambda\rightarrow F$ by
$$f_{g}(\theta)=(\mathcal{F}(\chi_{\{a_{\theta}\}}))(g), 
\quad
\mbox{for $\theta <\lambda$}.$$

\begin{claim}\label{c1}
For all $g\in G$, the net $(f_{g}(\theta))_{\theta <\lambda}$ is
convergent in $F$.
\end{claim}

\begin{proof}[Proof of Claim \ref{c1}]
Fix $g\in G$ and $\epsilon\in F^{>0}$. 
By the Cauchy condition for $\mathcal{F}$, there exists $\theta_{0}<\lambda$
such that $\forall \alpha,\beta\in[[F^{G}]]$ with 
$\alpha,\beta\geq\chi_{\{a_{\theta_{0}}\}}$, we have 
$|\mathcal{F}(\alpha)-\mathcal{F}(\beta)|<\epsilon\chi_{\{g\}}$. 
This shows that 
$|\mathcal{F}(\alpha)(g)-\mathcal{F}(\beta)(g)|<\epsilon$. 
Therefore $\forall\theta_{1},\theta_{2}\geq\theta_{0}$,
$$|f_{g}(\theta_{1})-f_{g}(\theta_{2})| = 
|(\mathcal{F}(\chi_{\{a_{\theta_{1}}\}}))(g) - 
(\mathcal{F}(\chi_{\{a_{\theta_{2}}\}}))(g)|<\epsilon,$$ 
since 
$\chi_{\{a_{\theta_{1}}\}},\chi_{\{a_{\theta_{2}}\}} 
\geq \chi_{\{a_{\theta_{0}}\}}$.
Hence the net $(f_{g}(\theta))_{\theta <\lambda}$ is Cauchy in $F$ and so
convergent there, since $F$ is Scott complete.
\end{proof}

Let $\gamma:G\rightarrow F$ be defined by 
$\gamma (g)=\lim(f_{g}(\theta))_{\theta <\lambda}$.

\begin{claim}\label{c2}
$(\forall \eta<\lambda)
(\exists \theta_{0}<\lambda)
(\forall \theta\geq\theta_{0})
(\forall g\in G^{<-a_{\eta}})
\mathcal{F}(\chi_{\{a_{\theta}\}})(g) = \gamma(g)$.
\end{claim}

\begin{proof}[Proof of Claim \ref{c2}]
For any $\epsilon=\chi_{\{-a_{\eta}\}}$, with $\eta <\lambda$, there
exists $\theta_{0}<\lambda$ such that 
$\forall\theta_{1},\theta_{2}<\lambda$ with 
$\theta_{1},\theta_{2}\geq\theta_{0}$, we have 
$$|\mathcal{F}(\chi_{\{a_{\theta_{1}}\}}) - 
\mathcal{F}(\chi_{\{a_{\theta_{2}}\}})|<\chi_{\{-a_{\eta}\}}.$$
This shows that for all $\theta_{1},\theta_{2}\geq\theta_{0}$ and
$g<-a_{\eta}$, we have 
$$(\mathcal{F}(\chi_{\{a_{\theta_{1}}\}}))(g) = 
(\mathcal{F}(\chi_{\{a_{\theta_{2}}\}}))(g)=\gamma (g).$$
Therefore, for all $\theta\geq\theta_{0}$ and $g\in G^{<-a_{\eta}}$, we
have $\mathcal{F}(\chi_{\{a_{\theta}\}})(g)=\gamma(g)$.
\end{proof}

\begin{claim}\label{c3}
$\gamma\in [[F^{G}]]$.
\end{claim}

\begin{proof}[Proof of Claim \ref{c3}]
It suffices to show that the support of $\gamma$ is well ordered. 
For all $g\in G$, there exists $\eta <\lambda$ such that
$g<-a_{\eta}$. Let $\theta_{0}$ be as in Claim \ref{c2}. 
As $g$ can not be the initial term of any infinite strictly decreasing
sequence in the support of $\mathcal{F}(\chi_{\{a_{\theta_{0}}\}})$,
Claim \ref{c2} shows that the same holds for the support of $\gamma$.
\end{proof}

\begin{claim}\label{c4}
The function $\mathcal{F}$ on $[[F^{G}]]$ tends to $\gamma$ at infinity.
\end{claim}

\begin{proof}[Proof of Claim \ref{c4}]
It is enough, by the Cauchy criterion for $\mathcal{F}$, to apply $\mathcal{F}$
on those $f$'s in $[[F^{G}]]$ that are of the form
$\chi_{\{a_{\theta}\}}$, for $\theta <\lambda$ and let $\theta$ tend to
$\lambda$. 
The result is then immediate from Claim \ref{c2}.
\end{proof}

The above claims give the result.
\end{proof}

\begin{lem}\label{ifthencf}
Suppose that $F$ is a Scott complete ordered field, $G$ is a 2-divisible
ordered abelian group, and $K$ is a dense ordered sub-field of $[[F^{G}]]$
which contains $F$. Assume furthermore that $K$ is closed under the
automorphism on $[[F^{G}]]$ sending $\chi_{\{g\}}$ to $\chi_{\{g+g\}}$ and
its inverse. 
Then $K$ satisfies $D_{\emptyset}S_{cf}COF$.
\end{lem}

\begin{proof}
If $K=[[F^{G}]]$, then the conclusion will trivially
hold since $[[F^{G}]]$ is Scott complete. Suppose $K\subsetneq
[[F^{G}]]$. Assume for the purpose of a contradiction that 
$\mathcal{F}$ is a $\emptyset$-definable divergent Cauchy function on $K$.
Fix a net $(k_{\alpha})_{\alpha<\lambda}$ of elements of $K$,
where $\lambda$ is cofinality of $[[F^{G}]]$, which is cofinal in
the latter. For any $f\in [[F^{G}]]$, let $\lambda (f)$ be the
least ordinal less than $\lambda$ such that $f<k_{\lambda (f)}$.
Consider the function $\tilde{\mathcal{F}}$ on $[[F^{G}]]$ with
$\tilde{\mathcal{F}}(f)=\mathcal{F}(k_{\lambda (f)})$. As a Cauchy
function on the Scott complete ordered field $[[F^{G}]]$, it will
converge to some $f_{0}\in [[F^{G}]]$. Clearly $f_{0}\not\in K$
and therefore $f_{0}\not\in F$. Let $\mathcal{A}$ be the ordered
field automorphism on $[[F^{G}]]$ described in the statement of
the lemma. The only fixed points of $[[F^{G}]]$ under $\mathcal{A}$
are elements of $F$ (since if there is a leading $F$-infinitely
large or a highest $F$-infinitely small term, then the result of
applying $\mathcal{A}$ to such elements will have a different leading
$F$-infinitely large, respectively highest $F$-infinitely small
term and so must be different). Therefore, 
$\mathcal{A}(f_{0}) \not = f_{0}$. 
Let 
$$C = 
\{f\in K : \mbox {$f$ is cofinally many times dominated in $K$ by 
values of $\mathcal{F}$}\}.$$ 
Now since in $[[F^{G}]]$, $\mathcal{F}$ tends to $f_{0}$, while 
$\mathcal{A}(\mathcal{F})$ approaches $\mathcal{A}(f_{0})$ there (as 
$\mathcal{A}$ is continuous) and also $K$ is dense in $[[F^{G}]]$, we get 
$C\not = \mathcal{A}(C)$. This contradicts $\emptyset$-definability of 
$\mathcal{F}$, since $\mathcal{A}$ restricted to $K$ is an (onto) 
automorphism.
\end{proof}

\begin{thm}\label{conditional}
\mbox{}
\begin{itemize}
\item[(i)]
$D_{p}S_{rc}COF\vdash D_{p}S_{cf}COF$, 
$D_{\emptyset}S_{rc}COF\vdash D_{\emptyset}S_{cf}COF$.
\item[(ii)]
$D_{\emptyset}S_{cf}COF\nvdash D_{p}S_{rc}COF$.
\end{itemize}
\end{thm}

\begin{proof}
(i). 
Let $F\vDash D_{p}S_{rc}COF$ (respectively
$F\vDash D_{\emptyset}S_{rc}COF$) and $f:F\rightarrow F$ be a
p-definable (respectively $\emptyset$-definable) Cauchy function.
Let $\psi(z)$ be the formula expressing, using $f$ as a shorthand,
$(\forall t)(\exists x\geq t)(z\leq f(x))$. We claim that the set
$C$ defined by $\psi$ in $F$ is a regular cut whose supremum is
the limit of $f$.

To see regularity of $C$, fix $\epsilon\in F^{>0}$. From the Cauchy
criterion for $f$, there exists $d\in F$ such that $\forall x,y\geq d$,
$|f(x)-f(y)|<\frac{\epsilon}{2}$. Let $z=f(d)-\frac{\epsilon}{2}$. Then
$F\vDash\psi (z)$ and $F\vDash\neg\psi (z+\epsilon)$.

Now let $\sup(C)=\alpha$. 
To show $\lim f(x)=\alpha$, as $x$
becomes arbitrarily large in $F$, let $\epsilon\in F^{>0}$ and $d$
be as before. 
By the definition of $\alpha$, there exists $z\in C$
with $z>\alpha-\frac{\epsilon}{2}$. 
From this, we get $\exists
x_{0}\geq d$ with $\alpha-\frac{\epsilon}{2}<z\leq f(x_{0})$. 
On
the other hand, for any $t\in F$, the element $x=\max\{t,d\}$
satisfies $f(x_{0})-\frac{\epsilon}{2}<f(x)$ and so
$f(x_{0})-\frac{\epsilon}{2}\leq\alpha$. 
Therefore,
$|f(x_{0})-\alpha|\leq\frac{\epsilon}{2}$. 
Hence, 
$$
\forall x\geq d,\ 
|f(x)-\alpha| \leq |f(x)-f(x_{0})|+|f(x_{0})-\alpha| 
< \frac{\epsilon}{2}+\frac{\epsilon}{2} 
= \epsilon.
$$

(ii). 
Let $R$ be the ordered ring $[[{\mathbb R}^{{\mathbb Q}^{\leq 0}}]]$ 
and $K$ its fraction field. We claim that $K$ is a witness
for the independence assertion at hand. By Lemmas 3.1 and 3.3 and
$[[{\mathbb R}^{\mathbb Q}]]$ being real closed, it suffices to
show that $K$ is a non real closed dense sub-field of
$[[{\mathbb R}^{\mathbb Q}]]$ closed under the automorphism
mentioned there and its inverse. The latter is immediate from the
same property for $R$.

As to $K$ being proper in its real closure, we use a result which has
recently been shown in \cite{P}. It states that in $R$, all elements $a+1$
where $a$ has a strictly increasing $\omega$-sequence support converging
to $0$, are prime. The series
$$t^{-1}+t^{-\frac{1}{2}}+t^{-\frac{1}{3}}+\cdots+1$$ 
is therefore an infinite prime in $R$. The square root of any (positive) 
infinite prime in $R$ (which has cofinal positive-integer powers there) 
witnesses that $K$ is not real closed.

To see that $K$ is dense in $[[{\mathbb R}^{\mathbb Q}]]$, notice that
$R$
approximates within $1$ any element of 
$[[{\mathbb R}^{\mathbb Q}]]$.
\end{proof}

It is interesting to note that, by Pitteloud's results again, the same 
infinite prime remains prime in all
$[[{\mathbb R}^{({\mathbb R}^{\alpha})^{\leq 0}}]]$ for any ordinal 
$\alpha$. 
The p-definable regular gap corresponding to its square root is therefore
never realized in the fraction fields of such ordered extension fields and
can not be traversed by $\emptyset$-definable functions over them.

\subsection*{Acknowledgement} 

The first author thanks Institute for Studies in Theoretical Physics and
Mathematics (IPM), Tehran, Iran, for partial support. This work will form
part of second author's PhD thesis under first author's supervision at
Tarbiat Modarres University, Tehran, Iran. The authors are grateful to
D. Pitteloud for kindly sending his paper to them before it was published.

\subsection*{Note added in proof}
The above proof for Proposition \ref{scottfieldsofgps} will be 
supplemented elsewhere by an appropriate argument to show that $[[F^{G}]]$
is always Scott complete, no matter whether $F$ is so or not (of course,
there $G$ will be a non-zero group). 
Lemma \ref{ifthencf} therefore holds without assuming that $F$ is Scott
complete. 
Theorem \ref{conditional}(ii) will also be improved to
$D_{\emptyset}S_{rc}COF\nvdash D_{p}S_{rc}COF$ (with the same witness $K$
as above). 
Consequently, either $D_{\emptyset}S_{cf}COF\nvdash^{(?)} D_{p}S_{cf}COF$
or $D_{p}S_{cf}COF\nvdash^{(?)} D_{p}S_{rc}COF$ (or both).
The parameter-free version of the latter remains open also.

\providecommand{\bysame}{\leavevmode\hbox to3em{\hrulefill}\thinspace}
\providecommand{\MR}{\relax\ifhmode\unskip\space\fi MR }
\providecommand{\MRhref}[2]{%
  \href{http://www.ams.org/mathscinet-getitem?mr=#1}{#2}
}
\providecommand{\href}[2]{#2}

\end{document}